\documentclass{amsart}

\usepackage{amsmath,amssymb}
\usepackage{graphicx}
\usepackage{tikz}
\usepackage{tikz-cd}

\newtheorem{theorem}{Theorem}
\newcommand{\bt}{\begin{theorem}}
\newcommand{\et}{\end{theorem}}
\newtheorem{lemma}{Lemma}
\newcommand{\bl}{\begin{lemma}}
\newcommand{\el}{\end{lemma}}
\newtheorem*{theoremNN}{Theorem}
\newcommand{\btNN}{\begin{theoremNN}}
\newcommand{\etNN}{\end{theoremNN}}
\newtheorem{problem}{Problem}
\newcommand{\bprob}{\begin{problem}}
\newcommand{\eprob}{\end{problem}}

\newcommand{\beq}{\begin{equation}}
\newcommand{\eeq}{\end{equation}}
\newcommand{\benum}{\begin{enumerate}}
\newcommand{\eenum}{\end{enumerate}}

\newcommand{\Z}{\ensuremath{\mathbf Z}}

\newcommand{\Fq}{\ensuremath{{\mathbf F}_q}}

\newcommand{\mcm}{\ensuremath{ \mathcal M}}

\newcommand{\R}{\ensuremath{\mathbf R}}

\newcommand{\mcs}{\ensuremath{ \mathcal S}}

\DeclareMathOperator{\qqand}{\qquad\text{and}\qquad}

\title[Estimates for linear forms]{Comparison estimates for linear forms 
in additive number theory}
\author{Melvyn B. Nathanson}
\address{Department of Mathematics\\Lehman College (CUNY)\\Bronx, NY 10468} \email{melvyn.nathanson@lehman.cuny.edu}

\subjclass[2010]{05A17, 11B13,  11B30, 11B75, 11P99.} 
\keywords{Sumsets, difference sets, linear forms, additive number theory.}
\thanks{Supported in part by a grant from the PSC-CUNY Research Award Program.}

\date{\today}

\begin{document}

\begin{abstract}
Let $R$ be a commutative ring $R$ with $1_R$ and with group of units  $R^{\times}$.   
Let $\Phi = \Phi(t_1,\ldots, t_h) = \sum_{i=1}^h \varphi_it_i$
be an $h$-ary linear form with nonzero coefficients $\varphi_1,\ldots, \varphi_h \in R$.  
Let $M$ be an $R$-module.  
For every subset $A$ of $M$,  the \emph{image of $A$ under $\Phi$} is 
\[
\Phi(A) = \{ \Phi(a_1,\ldots, a_h) : (a_1,\ldots, a_h) \in A^h \}.
\] 
For every  subset $I$ of $\{1,2,\ldots, h\}$, there is the 
\emph{subset sum} 
$
s_I = \sum_{i\in I} \varphi_i.  
$
Let   
$
\mcs(\Phi) = \{s_I:  \emptyset  \neq I \subseteq \{1,2,\ldots, h\} \}. 
$

\btNN      
Let $\Upsilon(t_1,\ldots, t_g) = \sum_{i=1}^g \upsilon_it_i$ and $\Phi(t_1,\ldots, t_h) = \sum_{i=1}^h \varphi_it_i$ be linear forms with nonzero coefficients in the ring $R$.  
If $\{0, 1\} \subseteq \mcs(\Upsilon)$ and $\mcs(\Phi) \subseteq R^{\times}$, 
then for every $\varepsilon > 0$  and $c > 1$ there exist a finite $R$-module $M$ 
with $|M| > c$ and a subset $A$ of $M$ such that 
$\Upsilon(A  \cup \{0\}) = M$ and 
$|\Phi(A)| < \varepsilon |M|$.
\etNN

\end{abstract}

\maketitle

\section{The problem}
In 1973, Haight~\cite{haig73} proved that for 
all positive integers $h$ and $\ell$ there exist a positive integer $m$ 
and a subset $A$ of $\Z/m\Z$ such that 
\[
A-A = \Z/m\Z
\]
but the $h$-fold sumset $hA$ omits $\ell$ consecutive congruence classes.  
Ruzsa~\cite{ruzs16}, refining Haight's method, recently proved that, 
for every positive integer $h$ and every $\varepsilon > 0$, there exist 
a positive integer $m$ and a subset $A$ of $\Z/m\Z$ such that 
\[
A-A = \Z/m\Z
\qqand |hA| < \varepsilon m.
\]
The difference set $A-A$ is the image of $A$ under the linear form 
$\Upsilon(t_1,t_2) = t_1-t_2$
and the $h$-fold sumset $hA$ is the image of $A$  under the linear form 
$\Phi(t_1,t_2, \ldots, t_h) = t_1 + t_2 + \cdots + t_h$.  
Equivalently, Ruzsa constructed a subset $A$ of the \Z-module $M = \Z/m\Z$ such that 
\[
\Upsilon(A) = M
\qqand |\Phi(A)| < \varepsilon |M|.
\]
This is a significant result in additive number theory.  
In this paper, we extend Ruzsa's theorem to a large class 
of pairs of linear forms $\Upsilon$ and $\Phi$.

Let $R$ be a commutative ring with multiplicative identity $1_R \neq 0$.
We denote the group of units in  $R$ by $R^{\times}$.  
Associated to every sequence $(\varphi_1,\ldots, \varphi_h)$ of nonzero elements of $R$ 
is the $h$-ary linear form  
\[
\Phi(t_1,\ldots, t_h) = \sum_{i=1}^h \varphi_it_i.
\]
For every  subset $I$ of $\{1,2,\ldots, h\}$, we define the 
\emph{subset sum} 
\[
s_I = \sum_{i\in I} \varphi_i.  
\]      
Note that $s_{\emptyset} = 0$ and  $s_{\{i\} } = \varphi_i$ for $i=1, \ldots, h$.  Let   
\[
\mcs(\Phi) = \{s_I:  \emptyset  \neq I \subseteq \{1,2,\ldots, h\} \} 
\]
be the set of all nonempty subset sums of the sequence of coefficients of $\Phi$.

For example, if $\Phi(t_1,\ldots, t_h) = t_1 + t_2 + \cdots + t_h$, then $S(\Phi) = \{1,2,\ldots, h\}$.
If $\Upsilon(t_1,t_2) = t_1 - t_2$, then $S(\Upsilon) = \{-1, 0, 1\}$.  

Let $M$ be an $R$-module.  
The linear form $\Phi$ induces the function $\Phi:M^h \rightarrow M$ defined by 
\[
\Phi(x_1,\ldots, x_h) = \sum_{i=1}^h \varphi_i x_i       
\]
for all $(x_1,\ldots, x_h) \in M^h$.     
For every subset $A$ of $M$, the \emph{image of $A$ under $\Phi$} is 
\[
\Phi(A) = \{ \Phi(a_1,\ldots, a_h) : (a_1,\ldots, a_h) \in A^h \}.
\]

In this paper we investigate the problem 
of  classifying the pairs of $R$-linear forms $(\Upsilon, \Phi)$ 
with the property that, for every $\varepsilon > 0$,  there exist a  
finite $R$-module $M$ and a subset $A$ of $M$ such that 
$\Upsilon(A) = M$ and $ |\Phi(A)| < \varepsilon |M|$.

A related problem for binary linear forms was previously investigated by 
Nathanson, O'Bryant, Orosz, Ruzsa, and Silva~\cite{nath-obry-oros-ruzs-silv07}.

\section{Results}
Let $R$ be a commutative ring, and let $M$ be an $R$-module.  
If $M \neq \{ 0\}$, then $|M| \geq 2$ and so 
$\lim_{k\rightarrow \infty} |M^k| = \lim_{k\rightarrow \infty} |M|^k = \infty$.
Thus,  if there is a nonzero finite $R$-module, then there 
are arbitrarily large finite $R$-modules.  
For $R = \Z$, the finite $\Z$-modules are the finite abelian groups.  For the finite field $\Fq$, 
the finite $\Fq$-modules are the finite-dimensional vector spaces over \Fq.
However, a nonzero vector space over an infinite field is infinite, 
and so an infinite field has no nonzero finite modules.  
In this paper we consider only commutative rings, like \Z\ and \Fq, 
for which nonzero finite modules exist.

\bt            \label{HR2:theorem:main}
Let $\Upsilon$ and $\Phi$ be linear forms with nonzero coefficients in the ring $R$.  
If 
\beq             \label{HR2:Upsilon-condition}
\{0, u\} \subseteq \mcs(\Upsilon) \quad \text{ for some } u \in R^{\times}
\eeq
and if 
\beq             \label{HR2:Phi-condition}
\mcs(\Phi) \subseteq R^{\times}
\eeq
then, for every $\varepsilon > 0$ and $c > 1$, there exist a finite $R$-module $M$ 
with $|M| > c$ and a subset $A$ of $M$ such that 
\beq             \label{HR2:main-conclusion}
\Upsilon(A \cup \{0\}) = M
\qqand 
|\Phi(A)| < \varepsilon |M|.  
\eeq
If $R$ is a finite field, the surjectivity 
condition $\Upsilon(A \cup \{0\}) = M$ 
can be replaced with $\Upsilon(A) = M$.

The construction of the $R$-module $M$ and the set $A$ depend 
only on the linear form $\Phi$ 
and not on the linear form $\Upsilon$.  
\et

\bt            \label{HR2:theorem:mainZ}
Let $\Upsilon$ and $\Phi$ be linear forms with nonzero integer coefficients.  
If 
\beq             \label{HR2:mainZ-condition}
0 \in \mcs(\Upsilon) \qqand 0 \notin \mcs(\Phi)
\eeq
then, for every $\varepsilon > 0$ and $c > 1$, there exist an integer $m > c$ 
and a subset $A$ of $\Z/m\Z$ such that 
\beq             \label{HR2:mainZ-conclusion}
\Upsilon(A ) = \Z/m\Z
\qqand 
|\Phi(A)| < \varepsilon m.  
\eeq
The construction of the integer $m$ and the set $A$ depend only 
on the linear form $\Phi$ and not on $\Upsilon$.   
\et

For example, let $R = \Z$ and let  $\Upsilon(t_1,t_2) = t_1 - t_2$ 
and $\Phi(t_1,\ldots, t_h) = \sum_{i=1}^h t_i$.  
We have  $0 \in \mcs(\Upsilon) = \{ -1, 0,1\}$ and $0 \notin \mcs(\Phi) = \{1,2,\ldots, h\}$,    
and so $\Upsilon$ and $\Phi$ satisfy the conditions of Theorem~\ref{HR2:theorem:mainZ}.  
This gives Ruzsa's result.  
Similarly, the linear forms $\Upsilon(t_1,t_2) = t_1 - t_2$ 
and $\Phi(t_1, t_2) = 2t_1 - t_2$ satisfy the conditions of Theorem~\ref{HR2:theorem:mainZ}.  
The linear forms $\Upsilon(t_1,t_2) = t_1 - t_2$ 
and $\Phi(t_1, \ldots, t_h) = 2t_1 +\sum_{i=2}^h  t_i$ also 
satisfy the conditions of Theorem~\ref{HR2:theorem:mainZ}.  
This answers a question in~\cite{nath16f}.  

We can extend Theorems~\ref{HR2:theorem:main} and~\ref{HR2:theorem:mainZ} 
to sets of three or more linear forms.   
For $k  \in \{ 1,\ldots, K\}$, let 
\[
\Phi_k(t_{1,k}, \ldots,t_{h_k,k} ) = \sum_{i=1}^{h_k} \varphi_{i,k} t_{i,k}
\]
be a linear form in $h_k$ variables with nonzero coefficients $ \varphi_{i,k}$ in  $R$.
The sum of these linear forms is the linear form $\chi =  \sum_{k=1}^K \Phi_k$:       
\[
\chi(t_{1,1},\ldots, t_{h_K,K}) 
= \sum_{k=1}^K  \Phi_k(t_{1,k}, \ldots,t_{h_k,k} ) 
= \sum_{k=1}^K  \sum_{i=1}^{h_k} \varphi_{i,k} t_{i,k}.  
\]
Thus, $\chi$ is a  linear form in $\sum_{k=1}^K  h_k$ variables with 
nonzero coefficients in $R$.

\bt                        \label{HR2:theorem:mainMany}
Let $\Upsilon_1,\ldots, \Upsilon_J$, and $\Phi_1,\ldots, \Phi_K$ be linear forms 
with nonzero coefficients in a ring $R$, 
and let $\chi =  \sum_{k=1}^K \Phi_k$.
If, for each $j \in \{1,\ldots, J \}$, there exists $u_j \in R^{\times}$ such that 
$\{0,u_j \} \subseteq \mcs(\Upsilon_j)$,  
and if $\mcs(\chi) \subseteq R^{\times}$, 
then, for every $\varepsilon > 0$ and $c > 1$, there exist a finite $R$-module $M$ 
with $|M| > c$ and a subset $A$ of $M$ such that 
\[
\Upsilon_j(A \cup \{ 0\}) = M
\]
for all $j=1,\ldots, J$, and
\[
|\Phi_k(A)| < \varepsilon |M|
\]
for all $k=1,\ldots, K$.
If $R$ is a finite field,  the surjectivity 
condition $\Upsilon_j(A \cup \{0\}) = M$ 
can be replaced with $\Upsilon_j(A) = M$ for all $j=1,\ldots, J$.
\et

\bt                          \label{HR2:theorem:mainManyZ}
Let $\Upsilon_1,\ldots, \Upsilon_J$, and $\Phi_1,\ldots, \Phi_K$ be linear forms 
with nonzero integer coefficients, 
and let $\chi =  \sum_{k=1}^K \Phi_k$.
If $0 \in \mcs(\Upsilon_j)$ for all $j \in \{1,\ldots, J\}$ and if $0 \notin \mcs(\chi)$, 
then, for every $\varepsilon > 0$ and $c >0$, there exist an integer $m > c$ 
and a subset $A$ of $\Z/m\Z$ such that 
\[
\Upsilon_j(A) = \Z/m\Z
\]
for all $j=1,\ldots, J$, and
\[
|\Phi_k(A)| < \varepsilon m
\]
for all $k=1,\ldots, K$.  
\et

The linear forms 
\[
\Upsilon(t_1) = t_1 \qqand \Phi(t_1, t_2) = t_1 + t_2       
\]
satisfy  
\[
0 \notin S(\Upsilon) = \{1\} \qqand S(\Phi) = \{1,2\}.
\]
Let $0 < \varepsilon < 1$.  For every $R$-module $M \neq \{ 0\}$, 
if $A$ is a subset of $M$ with 
$\Upsilon(A \cup \{0\} ) = M$, then $A = M$ or $M\setminus \{ 0\}$.
It follows that $\Phi(A) = M$ and $|\Phi(A)| > \varepsilon |M|$.   
Thus, the conclusion~\eqref{HR2:main-conclusion} 
of Theorem~\ref{HR2:theorem:main} does not necessarily 
apply to  linear forms $\Upsilon$ and $\Phi$ if $0 \notin S(\Upsilon)$.  

Let $R = \Z$, and consider the linear forms 
\[
\Upsilon(t_1,t_2) = t_1 + t_2 \qqand \Phi(t_1) = t_1    
\]
with  
\[
 0 \notin S(\Upsilon) = \{1,2\} \qqand S(\Phi) = \{1\}.
\]
For every positive integer $m$, there is a unique integer $d$ 
such that 
\[
\sqrt{m} \leq d < \sqrt{m}+1.
\]
Let 
\[
A_0 = \{ 0,1,2,\ldots, d-1\} \cup \{qd: q = 1 ,2,\ldots, d-1 \} \subseteq \Z.
\]
Every integer $n$ can be written uniquely in the form
$n = qd+r$, where $q \in \Z$ and $r \in  \{0,1,2,\ldots, d-1\}$.  
If $n \geq 0$, then $q \geq 0$.  
Because $d^2 \geq m$,  if $0 \leq n \leq m-1$, then $0 \leq q \leq d-1$.
It follows that 
\[
 \{ 0,1,2,\ldots, m-1\} \subseteq  A_0 + A_0 = \Upsilon(A_0).
\]
We also have $\Phi(A_0) = A_0$ and so
\[
|\Phi(A_0)| = 2d - 1 < 2\sqrt{m}+1.     
\]
In the finite $\Z$-module $M = \Z/m\Z$, let 
\[
A = \{ r+m\Z: r = 0,1,2,\ldots, d-1\} \cup \{qd + m\Z: q = 1 ,2,\ldots, d-1 \} \subseteq \Z/m\Z.
\]
We have
\[
\Upsilon(A) = \Z/m\Z \qqand \Phi(A) = A.
\]
Let $\varepsilon > 0$.  For all sufficiently large $m$, we have 
\[
|\Phi(A)| \leq 2d-1 <  2\sqrt{m}+1 < \varepsilon m.
\]
Thus, conditions~\eqref{HR2:Upsilon-condition} and~\eqref{HR2:Phi-condition} 
of Theorem~\ref{HR2:theorem:main} and 
condition~\eqref{HR2:mainZ-condition} of Theorem~\ref{HR2:theorem:mainZ} 
are not necessary, and the conclusions of Theorem~\ref{HR2:theorem:main} 
and~\ref{HR2:theorem:mainZ} may also 
apply to a pair of linear forms $(\Upsilon, \Phi)$ with $0 \notin S(\Upsilon)$.

\section{Surjectivity conditions}
The following lemma proves the surjectivity result for the linear form $\Upsilon$ 
 in Theorem~\ref{HR2:theorem:main}
and the linear forms $\Upsilon_j$ in Theorem~\ref{HR2:theorem:mainMany}.

\bl                       \label{HR2:lemma:surjective}
Let $M$ be an $R$-module.  
For every function $f:M \rightarrow M$, 
where $f$ is not necessarily an $R$-module homomorphism, let 
\[
A = A(M,f) =  \{ f(x) : x \in M \} \cup \{ f(x) +  x: x \in M \}.  
\]
Let 
\[
\Upsilon(t_1,\ldots, t_g) = \sum_{i=1}^g \upsilon_i t_i
\]
be a $g$-ary linear form with nonzero coefficients $\upsilon_1,\ldots, \upsilon_g$ in $R$ such that 
there exist nonempty subsets $I$ and $J$ of $\{1,2,\ldots, g\}$ 
with $s_I = 0$ 
and $s_J \in R^{\times}$.  
If $I \cup J \neq \{1,2,\ldots, g \}$, then  
\[
\Upsilon(A \cup \{ 0 \}) = M.
\]
If $I \cup J = \{1,2,\ldots, g\}$ or if $R$ is a  field, then 
\[
\Upsilon(A) = M.
\]
\el

\begin{proof}
The subset sum $s_J$ is a unit in $R$, and so we can define the linear form 
\[
\Upsilon'(t_1,\ldots, t_g) = s_J^{-1}\Upsilon(t_1,\ldots, t_g) 
= \sum_{i=1}^g s_J^{-1} \upsilon_i t_i  = \sum_{i=1}^g  \upsilon'_i t_i  
\]
where
$
 \upsilon'_i = s_J^{-1} \upsilon_i 
$
for $i \in \{1,\ldots, g\}$.
We have the subset sums 
\[
s'_I = \sum_{i\in I} \upsilon'_i =  s_J^{-1}\sum_{i\in I} \upsilon_i 
=  s_J^{-1}s_I = s_J^{-1}\cdot 0 = 0
\]
and
\[
s'_J = \sum_{j \in J} \upsilon'_j =  s_J^{-1}\sum_{j \in J} \upsilon_j 
=  s_J^{-1}s_J= 1.
\]
Because $\Upsilon'(A) = M$ if and only if $\Upsilon(A) = M$, we can assume that $s_J = 1$.  

For  $x \in M$, let  
\[
x_i = \begin{cases}
f(x) & \text{ if $i \in I\setminus J$}  \\
f(x) +  x& \text{ if $i \in I \cap J$}  \\
 x & \text{ if $i \in J \setminus I$}  \\
0 & \text{ if $i \in \{1,\ldots, g\} \setminus (I \cup J)$.}
\end{cases}
\]
We obtain  
\begin{align*}
\Upsilon(x_1,\ldots, x_g) 
& =  \sum_{i\in   I\setminus J} \upsilon_i x_i +  \sum_{i\in   I \cap  J} \upsilon_i x_i +  \sum_{i\in  J \setminus I} \upsilon_i x_i  +  \sum_{i \notin   I \cup  J} \upsilon_i x_i \\ 
& = \sum_{i\in I\setminus J} \upsilon_i  f(x)+  \sum_{i\in   I \cap  J} \upsilon_i (f(x)+  x) 
+  \sum_{i\in  J \setminus I} \upsilon_i   x \\
& = \sum_{i\in I} \upsilon_i  f(x) +  \sum_{i\in  J } \upsilon_i  x \\
& = s_I f(x) + s_J x  \\
& = 0 \cdot f(x) + 1 \cdot x  \\
& = x.
\end{align*}
It follows that $\Upsilon(A \cup \{0\}) = M$.  If $I \cup J = \{1,\ldots, g\}$, 
then $\Upsilon(A) = M$. 

Let $R$ be a field.   We have $0 \neq \upsilon_1 = s_{ \{1\} } \in R^{\times}$. 
If $\upsilon^* = \sum_{i=1}^g \upsilon_i = 0$, then we can choose 
$I = \{1,\ldots, g\}$ and $J = \{ 1\}$.  
If $\upsilon^* \neq 0$, then $\upsilon^* \in R^{\times}$ and 
we can choose $I$ such that $s_I = 0$ and $J = \{1, \ldots, g\}$.  In both cases, 
$I \cup J = \{1, \ldots, g\}$ and $\Upsilon(A) = M$.  
This completes the proof.  
\end{proof}

The following result gives the surjectivity parts of 
Theorems~\ref{HR2:theorem:mainZ}
and~\ref{HR2:theorem:mainManyZ}.

\bl                       \label{HR2:lemma:surjectiveZ} 
For every function $f: \Z/m\Z  \rightarrow \Z/m\Z $, 
where $f$ is not necessarily a group homomorphism, 
let 
\[
A = A(\Z/m\Z,f) =  \{ f(x) : x \in \Z/m\Z \} \cup \{ f(x) + x: x \in \Z/m\Z \}.  
\]
Let 
\[
\Upsilon(t_1,\ldots, t_g) = \sum_{i=1}^g \upsilon_i t_i
\]
be a $g$-ary linear form with nonzero  integer coefficients  such that 
$0 \in \mcs(\Upsilon)$.  
If $\gcd(s,m) = 1$ for all $s \in \mcs(\Upsilon)\setminus \{ 0\}$, 
then 
\[
\Upsilon(A) = \Z/m\Z.
\]
\el

\begin{proof}
Let
\[
\upsilon^* = \sum_{i=1}^g \upsilon_i \in \Z.
\]
If $\upsilon^* = 0$, let $I = \{1,2,\ldots, g\}$ and let $J = \{1\}$.  
Then $s_I = \upsilon^* = 0$ and $s_J = \upsilon_1 \neq 0$.  
If $\upsilon^* \neq 0$, let $J = \{1,2,\ldots, g\}$ and let $I$ be a subset of $\{1, 2,\ldots, g\}$ 
such that $s_I = 0$.  
Then $s_I =0$ and $s_J = \upsilon^* \neq 0$.  

In both cases, we have  $s_I = 0$, $s_J \neq 0$, 
and $I \cup J = \{1,2,\ldots, g\}$.
Because $\gcd(s_J,m) = 1$, it follows that $s_J$ is a unit in the ring $R = \Z/m\Z$.
With $M = \Z/m\Z$, we have  $\Upsilon(A) = \Z/m\Z$ by Lemma~\ref{HR2:lemma:surjective}.  
This completes the proof.  
\end{proof}

\section{Admissible pairs of functions}

Consider the  linear form 
\[
\Phi(t_1,\ldots, t_h) = \sum_{i=1}^h \varphi_it_i
\]
with nonzero coefficients  $\varphi_1, \ldots, \varphi_h \in R$.  
Throughout this section we assume that 
\[
0 \notin S(\Phi).
\]
Let $h$ and $\ell$ be positive integers with $\ell \leq h$.
Let 
\beq            \label{HR2:part1}
\{1,2,\ldots, h \}  =  I_1 \cup I_2 \cup \cdots \cup I_{\ell}  
\eeq
is a partition of $\{1,2,\ldots, h\}$ into $\ell$ pairwise disjoint nonempty sets.   
For $j=1,\ldots, \ell$, let  
\beq            \label{HR2:part2}
I_{j} = I_{j,0} \cup I_{j,1}  \qqand  I_{j,0} \cap I_{j,1}  = \emptyset. 
\eeq
We do not assume that both sets $ I_{j,0}$ and $I_{j,1}$ are nonempty.

For $j=1,\ldots, \ell$, we have the  subset sums 
\begin{align*}
s_{I_{j,0}} & = \sum_{i\in I_{j,0}} \varphi_i \\
s_{I_{j,1}} & = \sum_{i\in I_{j,1}} \varphi_i \\
s_{I_j} & = \sum_{i\in I_j} \varphi_i  = s_{I_{j,0}} + s_{I_{j,1}}.        
\end{align*}

Let $M$ be an $R$-module.
A pair of functions $(\alpha, \beta)$, where 
\[
\alpha:M \rightarrow S(\Phi) \cup \{ 0\} 
\]
and 
\[
\beta:M \rightarrow S(\Phi) \cup \{ 0\}, 
\]
is \emph{admissible} if, for positive integers $\ell \leq h$,  
there is a set $\{y_1,\ldots, y_{\ell}\}$ of $\ell$ distinct elements of $M$ and a partition
of the set $\{1,\ldots, h\}$ of the form~\eqref{HR2:part1} and~\eqref{HR2:part2} 
such that, for all $x \in M$,  
\beq                 \label{HR2:alpha}
\alpha(x) = 
\begin{cases}
s_{I_{j,0}} & \text{if $x = y_j$ for some $j \in \{1,\ldots, \ell \}$ } \\ 
0               & \text{if $x \notin  \{y_1,\ldots, y_{\ell}\}$} 
\end{cases}
\eeq
and
\beq                 \label{HR2:beta}
\beta(x) = 
\begin{cases}
s_{I_{j,1}} & \text{if $x = y_j$ for some $j \in \{1,\ldots, \ell \}$ } \\ 
0               & \text{if $x \notin  \{y_1,\ldots, y_{\ell}\}$.} 
\end{cases}
\eeq

Let $M$ be an $R$-module, let $f:M \rightarrow M$ be a function, and let 
\beq          \label{HR2:AMf}
A = A(M,f) =  \{ f(x) : x \in M \} \cup \{ f(x) + x: x \in M \}.   
\eeq
If   
\[
w \in \Phi(A)
\]
then there exist sequences  
\[
(x_1,\ldots, x_h) \in M^h \qqand 
(\lambda_1, \ldots, \lambda_h) \in \{0,1\}^h
\]
such that 
\beq          \label{HR2:wPhi}
w  = \Phi \left(  f(x_1) + \lambda_1 x_1, \ldots,  f(x_h) + \lambda_h  x_h \right).  
\eeq
This representation of $w$ by the linear form $\Phi$ has  \emph{level} $\ell$ if 
\[
\ell = | \{x_1,\ldots, x_h\}|.   
\]
Note that an element $w \in \Phi(A)$ can have many representations.
For example, let $R = \Z$ and $M = \Z/10\Z$.  If $f(x) = 0$, then $A = \Z/10\Z$.  
Choosing $\Phi = t_1 + t_2 + t_3$, we have
 \begin{align*}
 9 + 10\Z & = \Phi(3+ 10\Z ,3+ 10\Z ,3+ 10\Z ) \\
&  = \Phi(1+ 10\Z ,4+ 10\Z ,4+ 10\Z ) \\ 
& = \Phi(2+ 10\Z ,3+ 10\Z ,4+ 10\Z ).
 \end{align*}
These are representations of $9+ 10\Z$  of levels 1,2, and 3, respectively.  

Let $w \in \Phi(A)$ have the representation~\eqref{HR2:wPhi} of level $\ell$, 
and let 
\[
\{x_1,\ldots, x_h \} = \{ y_1, \ldots, y_{\ell} \}.
\]
For $j=1,\ldots, \ell$, let 
\begin{align*}
I_{j} & = \{ i \in \{1,\ldots, h\}: x_i = y_j \}  \\
I_{j,0} & = \{ i \in I_j:  \lambda_i = 0\} \\
I_{j,1} & = \{ i \in I_j : \lambda_i = 1\}.  
\end{align*}
We obtain a partition of the set $\{1,\ldots, h\}$ 
of the form~\eqref{HR2:part1} and~\eqref{HR2:part2}, 
with the associated subset sums $s_{I_{j,0}}$, $s_{I_{j,1}}$, 
and  $s_{I_j}$.  
It follows that 
\begin{align*}
w & = \Phi \left(  f(x_1) + \lambda_1 x_1, \ldots,  f(x_h) + \lambda_h  x_h \right) \\
& = \sum_{i=1}^h \varphi_i  \left(  f(x_i) + \lambda_i x_i \right)   \\
& = \sum_{j=1}^{\ell}  
\sum_{i \in I_j} \varphi_i \left(  f(y_j) + \lambda_j y_j  \right)  \\
& = \sum_{j=1}^{\ell} \left( 
\sum_{i \in I_{j,0}} \varphi_i  f(y_j)  
+ \sum_{i \in I_{j,1}} \varphi_i \left(  f(y_j) + y_j  \right) 
\right)  \\
 & = \sum_{j=1}^{\ell} \left( s_{I_{j,0}} f(y_j) +  s_{I_{j,1}} \left(  f(y_j) + y_j  \right)   \right) \\
& = \sum_{j=1}^{\ell}   \left( \alpha(y_j) f(y_j) +  \beta(y_j) \left(  f(y_j) + y_j  \right)    \right) \\
& = \sum_{x \in M} \left( \alpha(x) f(x) + \beta(x) (f(x)+x) \right)
\end{align*}
where the functions $\alpha$ and $\beta$ are defined by~\eqref{HR2:alpha} 
and~\eqref{HR2:beta}, and the pair of functions $(\alpha, \beta)$ is admissible.  

Conversely, let $(\alpha, \beta)$ 
be an admissible pair of functions associated with a partition 
of $\{1,\ldots, h\}$ of the form~\eqref{HR2:part1} and~\eqref{HR2:part2} 
and a subset $\{y_1,\ldots, y_{\ell}\}$ of $M$ of cardinality $\ell$.
Define 
\[
x_i = y_j \qquad \text{if $i \in I_j$}
\]
and
\[
\lambda_i = 
\begin{cases}
0 & \text{if $i \in I_{j,0}$}  \\ 
1 & \text{if $i \in I_{j,1}$.} 
\end{cases}
\]
For every function $f:M \rightarrow M$ we have 
\begin{align*}
w & = \sum_{x\in M} \left( \alpha(x) f(x) + \beta(x)(f(x) + x)  \right)  \\
& = \sum_{j=1}^{\ell} \left( \alpha(y_j)f(y_j)  + \beta(y_j) (f(y_j)  +  y_j) \right) \\
& = \sum_{j=1}^{\ell} \left( s_{I_{j,0}}f(y_j)  + s_{I_{j,1}} (f(y_j)  +  y_j) \right) \\
& = \sum_{j=1}^{\ell} \left( \sum_{i \in I_{j,0}} \varphi_i  f(x_i)  
+  \sum_{i \in I_{j,1}} \varphi_i  (f(x_i)  +  x_i) \right) \\
& =  \sum_{j=1}^{\ell}  \sum_{i \in I_j } \varphi_i  ( f(x_i)  +  \lambda_i x_i) \\
& =  \sum_{i=1}^h \varphi_i  ( f(x_i)  +  \lambda_i x_i) \\
& \in \Phi(A(M,f))
\end{align*}
where $A(M,f)$ is the subset of $M$ defined by~\eqref{HR2:AMf}.  
Thus, every admissible pair of functions $(\alpha, \beta)$  
with support $\{y_1,\ldots, y_{\ell} \}$ in $M$ determines a representation of level $\ell$ 
of an element $w \in \Phi(A(M,f))$.

\bl                      \label{HR2:lemma:commute}
Let $M'_0, M'_1, \ldots, M'_n$ be finite $R$-modules, and let 
\[
M = \bigoplus_{i=0}^n M'_i = \{ x=(x_0,x_1,\ldots, x_n):x_i \in M'_i \text{ for } i=0,1,\ldots, n \}.
\]  
Define the $R$-linear projection $\pi_0:M \rightarrow M'_0$ by $\pi_0(x) = x_0$.  

Let $f: M \rightarrow M$ and $f_0:M'_0 \rightarrow M'_0$ be functions such that 
$\pi_0(f(x)) = f_0(\pi_0(x))$  for all $x \in M$.  
Equivalently, the diagram 
\begin{center}
\begin{tikzcd}
M \arrow{r}{f}   \arrow{d}[swap]{\pi_0}  & M   \arrow{d}{\pi_0} \\
M'_0 \arrow{r}{f_0}  & M'_0 
\end{tikzcd}
\end{center}
commutes. 
Let 
\[
A = A(M,f) =  \{ f(x) : x \in M \} \cup \{ f(x) + x: x \in M \} 
\]
and
\[
A_0 = A(M'_0,f_0) =  \{ f_0(x_0) : x_0 \in M'_0 \} \cup \{ f_0(x_0) + x_0 : x_0 \in M'_0 \}.   
\]

Let $\Phi:M^h \rightarrow M$ be the  $h$-ary linear form defined by 
\[
\Phi(t_1,\ldots, t_h) = \sum_{i=1}^h \varphi_it_i
\]
with nonzero coefficients $\varphi_1,\ldots, \varphi_h \in R$.  
Then
\benum
\item[(i)]
\[
\pi_0(\Phi(A)) = \Phi(A_0).
\]
\item[(ii)]
Let $w \in \Phi(A)$ be represented by the admissible pair of functions 
$(\alpha,\beta)$  on $M$.  
Define functions $\alpha^*$ and $\beta^*$ from $M'_0$ to $S(\Phi) \cup \{ 0\}$
by   
\[
\alpha^*(x_0) = \sum_{ \substack{x\in M \\ \pi_0(x) = x_0 }} \alpha(x) 
\]
and
\[
\beta^*(x_0) = \sum_{ \substack{x\in M \\ \pi_0(x) = x_0 }} \beta(x)   
\] 
for all $x_0 \in M'_0$.
The pair $(\alpha^*, \beta^*)$ is admissible and represents $\pi_0(w)$. 
\item[(iii)]
If $(\alpha,\beta)$ has level $\ell$ and support $\{y_1,\ldots, y_{\ell}\}$, 
then $(\alpha^*, \beta^*)$ has level $\ell^* \in \{1, \ldots, \ell \}$ and support 
$\{z_1,\ldots, z_{\ell^*} \} = \{  \pi_0(y_1),\ldots, \pi_0(y_{\ell})  \}$.
\eenum
\el

\begin{proof} 
Let $w \in \Phi(A)$.  There exist $(x_1,\ldots, x_h) \in M^h$ 
and $(\lambda_1,\ldots, \lambda_h) \in \{0,1\}^h$ such that 
\[
w = \Phi\left( f(x_1)+\lambda_1 x_1, \ldots,   f(x_h)+\lambda_h x_h\right)     
= \sum_{i=1}^h \varphi_i (f(x_i)+\lambda_i x_i).
\]
For all  $i=1,\ldots, h$ and $j = 0,1,\ldots, n$, there exist elements $x_{i,j} \in M'_j$ 
such that 
\[
x_i = (x_{i,0}, x_{i,1}, \ldots, x_{i,n}) \in \bigoplus_{j=0}^n M'_j = M
\]
and
\begin{align*}
\pi_0(w) & = \pi_0\left(     \sum_{i=1}^h \varphi_i (f(x_i)+\lambda_i x_i)   \right) \\
 & =   \sum_{i=1}^h \varphi_i \left(    \pi_0 f(x_i)+\lambda_i \pi_0(x_i)   \right) \\
 & =   \sum_{i=1}^h \varphi_i \left(  f_0  \pi_0 (x_i)+\lambda_i x_{i,0}  \right) \\
  & =   \sum_{i=1}^h \varphi_i \left(  f_0 (x_{i,0})+\lambda_i x_{i,0}  \right) 
\in \Phi(A_0).
\end{align*}
Therefore, $\pi_0(\Phi(A)) \subseteq \Phi(A_0)$.  

Conversely, let $w_0 \in \Phi(A_0)$.  
There exist $(x_{0,1},\ldots, x_{0,h} )\in (M'_0)^h$ 
and $(\lambda_1,\ldots, \lambda_h) \in \{0,1\}^h$ such that 
\[
w_0 = \Phi\left( f_0(x_{0,1})+\lambda_1 x_{0,1}, \ldots,   f_0(x_{0,h})+\lambda_h x_{0,h} \right)     
= \sum_{i=1}^h \varphi_i (f_0(x_{i,0})+\lambda_i x_{i,0}).
\]
For $i=1,\ldots, h$, let 
\[
x_i =  (x_{i,0}, 0, \ldots, 0) \in \bigoplus_{j=0}^n M'_j = M
\]
and 
\[
w = \Phi\left( f(x_1)+\lambda_1 x_1, \ldots,   f(x_h)+\lambda_h x_h\right)     
= \sum_{i=1}^h \varphi_i (f(x_i)+\lambda_i x_i).
\]
It follows as above that 
\[
 w_0 =   \sum_{i=1}^h \varphi_i \left(  f_0 (x_{0,i})+\lambda_i x_{0,i}  \right) 
 = \pi_0(w) \in \pi_0(\Phi(A)) 
\]
and so $\Phi(A_0) \subseteq \pi_0(\Phi(A))$.
This proves~(i).

If $w \in \Phi(A)$ is represented by the admissible pair $(\alpha, \beta)$, then 
\[
w =  \sum_{x \in M} \left(  \alpha(x) f(x) 
+ \beta(x) (f(x)  + x) \right) 
\]
and so
\begin{align*}
\pi_0(w) & = \pi_0 \left( \sum_{x \in M} \left(  \alpha(x) f(x) 
+ \beta(x) (f(x)  + x) \right)  \right) \\
& =  \sum_{x \in M} \left(  \alpha(x) \pi_0 f(x) + \beta(x)  \pi_0 (f(x)  + x) \right) \\
& =  \sum_{x \in M} \left(  \alpha(x) f_0 \pi_0(x) + \beta(x) (  f_0\pi_0(x) + \pi_0(x)) \right) \\
& =  \sum_{x_0 \in M'_0} \left(  \left(\sum_{ \substack{x\in M \\ \pi_0(x) = x_0 }} 
\alpha(x)\right)  f_0(x_0) +  \left(\sum_{ \substack{x\in M \\ \pi_0(x) = x_0 }} 
\beta(x)\right)  (  f_0(x_0) + x_0) \right) \\ 
& =  \sum_{x_0 \in M'_0} \left(  \alpha^*(x_0)  f_0(x_0) +  \beta^*(x_0)  (  f_0(x_0) + x_0) \right). 
\end{align*}

Associated with the admissible pair $(\alpha, \beta)$ are partitions
\[
\{ 1,\ldots, h\} = I_1 \cup \cdots \cup I_{\ell} \qqand I_j = I_{j,0} \cup I_{j,1}
\]
and a set $\{y_1,\ldots, y_{\ell} \}$ of  $\ell$ distinct elements of $M$.  
Let 
\[
 \{  \pi_0(y_1),\ldots, \pi_0(y_{\ell})  \}=  \left\{z_1,\ldots, z_{\ell^*} \right\}  \subseteq M'_0
\]
 where $\ell^* \in \{ 1,2,\ldots, \ell \}$ 
and the elements $z_1,\ldots, z_{\ell^*}$ are distinct.

For $k \in \{1,\ldots, \ell^*\}$, we define  
\[
Q_k = \left\{ j \in \{ 1,\ldots, \ell \} : \pi_0(y_j) = z_k \right\}
\]
and 
\[
I^*_k = \bigcup_{j\in Q_k} I_j = \{  i \in \{1,\ldots, h\} :  \pi_0(x_i) = z_k   \}.
\]
It follows that 
\[
\{1,\ldots, \ell \} = Q_1 \cup \cdots \cup Q_{\ell^*}  
\]
and 
\[
 \{1,\ldots, h\} = I^*_1 \cup \cdots \cup I^*_{\ell} 
\]
are partitions into  pairwise disjoint nonempty sets.  
Let  
\[
I^*_{k,0} = \bigcup_{j \in Q_k} I_{j,0} \qqand I^*_{k,1} = \bigcup_{j \in Q_k} I_{j,1}.      
\]
We have  
\[
I^*_k = I^*_{k,0} \cup I^*_{k,1} \qqand  I^*_{k,0} \cap I^*_{k,1} = \emptyset.
\]
The partition $ \{1,\ldots, h\} = I^*_1 \cup \cdots \cup I^*_{\ell} $ 
and the set  $\left\{z_1,\ldots, z_{\ell^*} \right\}$ determine an admissible pair 
of functions $(\hat{\alpha}, \hat{ \beta})$ from $M'_0$ into $\mcs(\Phi) \cup \{ 0\}$ 
as follows:  
For all $x_0 \in M'_0$,  
\[
\hat{\alpha}(x_0) = 
\begin{cases}
s_{I^*_{k,0}}   & \text{if $x_0 = z_k$ for some $k \in \{1,\ldots, \ell^*\}$ } \\
0 &  \text{if $x_0 \notin \{ z_1, \ldots, z_{\ell^*} \}$.} 
\end{cases}
\]
and 
\[
\hat{\beta}(x_0) = 
\begin{cases}
s_{I^*_{k,1}}   & \text{if $x_0 = z_k$ for some $k \in \{1,\ldots, \ell^*\}$ } \\
0 &  \text{if $x_0 \notin \{ z_1, \ldots, z_{\ell^*} \}$.} 
\end{cases}
\]
For $k \in \{1,\ldots, \ell^*\}$, we have
\begin{align*}
\hat{\alpha}(z_k) & =  s_{I^*_{k,0}}  = \sum_{i \in I^*_{k,0}}  \varphi_i  
 =  \sum_{i \in \bigcup_{j \in Q_k} I_{j,0}   }  \varphi_i  \\
&  = \sum_{j \in Q_k} \sum_{i \in \ I_{j,0}   }  \varphi_i  
 =  \sum_{ \substack{ j=1 \\ \pi_0(y_j) = z_k    }   }^{\ell}  s_{I_{j,0}} \\
& =  \sum_{ \substack{ j=1 \\ \pi_0(y_j) = z_k    }   }^{\ell}  \alpha(y_j)  
 =  \sum_{ \substack{ x \in M \\ \pi_0(x) = z_k   } }  \alpha(x)  \\
& = \alpha^*(z_k)
\end{align*}
and so
\[
\hat{\alpha}(x_0) = \alpha^*(x_0). 
\]
Similarly, 
\[
\hat{\beta}(x_0) = \beta^*(x_0)
\]
for all $x_0 \in M'_0$.  
Therefore, $\pi_0(w) \in \Phi(A_0)$, 
and $(\alpha^*, \beta^*)$ is an admissible pair of functions on $M'_0$ 
of level $\ell^*$ and with support $\{z_1,\ldots, z_{\ell^*} \}$.   
This completes the proof.
\end{proof}

\section{Proof of Theorem~\ref{HR2:theorem:main}: The initial step}
We must prove that there exist $M$ and $A$ such that $|\Phi(A)| < \varepsilon |M|$,      
where  
\[
\Phi = \Phi(t_1,\ldots, t_h) = \sum_{i=1}^h \varphi_i t_i
\]
is an $h$-ary linear form with nonzero coefficients $\varphi_i \in R$ and 
subset sums $S(\Phi) \subseteq R^{\times}$.  

Let $\ell \in \{ 1,2,\ldots, h \}$.   
For every subset $A$ of an $R$-module $M$, we define 
\[
\Phi^{(\ell)}(A) = \{w \in \Phi(A): \text{$w$ has a representation of level at most $\ell$}\}.
\]
Because $\Phi$ is a function of $h$ variables, we have  $\Phi^{(h)}(A) = \Phi(A)$.  

Let $\varepsilon > 0$, and choose $\varepsilon_1, \ldots, \varepsilon_h$ such that 
\[
0 < \varepsilon_1 < \varepsilon_2 < \cdots < \varepsilon_h < \varepsilon.
\]
For all $\ell \in \{ 1,2,\ldots, h \}$, we shall construct  
a finite $R$-module $M_{\ell}$  
and a function $f_{\ell}:M_{\ell} \rightarrow M_{\ell}$ such that the set 
\[
A_{\ell} = A(M_{\ell},f_{\ell}) 
= \left\{ f_{\ell}(x): x \in M_{\ell} \right\} \cup \left\{ f_{\ell}(x) + x: x \in M_{\ell} \right\}
\]
satisfies 
\[
|\Phi^{(\ell)}(A_{\ell})| < \varepsilon_{\ell} |M_{\ell}|.
\]
Choosing $M = M_h$ and $A = A_h$, we obtain
\[
|\Phi(A)| = |\Phi^{(h)}(A_h)| < \varepsilon_h |M_h| < \varepsilon |M|.
\]
From Lemma~\ref{HR2:lemma:surjective}, we have 
$\Upsilon(A_{\ell} \cup \{ 0\}) = M_{\ell}$ for all $\ell = 1, \ldots, h$, 
and so $\Upsilon(A \cup \{ 0\}) = M$.   

The proof is by induction on $\ell$.
We begin with the case $\ell = 1$ and the construction of the module $M_1$ 
and the function $f_1:M_1 \rightarrow M_1$.

For every finite $R$-module $M$, function $f:M \rightarrow M$, 
and subset $A = A(M,f)$ of $M$, we have 
$w \in \Phi^{(1)}(A)$ if and only if there exist $x \in M$ and 
$(\lambda_1, \ldots, \lambda_h)  \in \{0,1\}^h$ such that 
\begin{align*}
w & = \sum_{i=1}^h \varphi_i ( f(x) +\lambda_i x) \\
&  = \left( \sum_{i=1}^h \varphi_i \right) f(x) +  \left( \sum_{i=1}^h  \lambda_i \varphi_i \right) x\\
& = \varphi^* f(x) + s_I x 
\end{align*}
where  
\[
\varphi^* = \sum_{i=1}^h \varphi_i 
\]
\[
I = \{i\in \{1,2,\ldots, h\}: \lambda_i = 1\}  
\]
and 
\[
s_I \in  \mcs(\Phi) \cup \{s_{\emptyset} \} = \mcs(\Phi) \cup \{ 0\} .
\]

For all $s \in  \mcs(\Phi) \cup \{ 0\} $, 
let $M'_{s}$ be a finite $R$-module such that 
\beq               \label{HR2:cardinalityM1} 
|M'_{s}| >  \max\left(  \frac{|\mcs(\Phi)| + 1}{\varepsilon_1} ,c \right)
\eeq
and let
\[
M_1 = \bigoplus_{{s} \in \mcs(\Phi) \cup \{ 0\} } M'_{s}. 
\]
Note that the construction of the finite module $M_1$ depends only on  
the set of subset sums of  the linear form $\Phi$, 
and not on the linear form $\Upsilon$.

If $x \in M_1$, then $x = (x_{s})_{s \in \mcs(\Phi) \cup \{ 0\} }$, where  
$x_{s} \in M'_{s}$.   For all $s \in \mcs(\Phi) \cup \{ 0\} $,  
we define the projection $\pi_{s}:M_1 \rightarrow M'_{s}$ by 
\[ 
\pi_{s}(x) = x_{s}
\]
and we define the function $g_{s}:M_1 \rightarrow M'_{s}$ by
\beq                     \label{HR2:gs} 
g_s(x)  = -\frac{{s}}{\varphi^*} \pi_{s}(x) = -\frac{{s}}{\varphi^*} \   x_{{s}}.
\eeq
This is possible because $\varphi^* \in \mcs(\Phi) \subseteq R^{\times}$.

Consider the function $f_1:M_1 \rightarrow M_1$ defined by
\beq                     \label{HR2:f1} 
f_1(x) = (g_s(x))_{s \in  \mcs(\Phi) \cup \{ 0\}} 
= \left(  -\frac{s}{\varphi^*} \   x_{s} \right)_{s \in \mcs(\Phi) \cup \{ 0\} } 
\eeq
and let 
\[
A = A(M_1,f_1) = \{ f_1(x):x\in M_1 \} \cup  \{ f_1(x) + x :x\in M_1 \}.  
\]

If $w \in \Phi^{(1)}(A)$, then there exist $x \in M_1$ 
and $s_I \in \mcs(\Phi) \cup \{ 0\} $ such that 
\[
w = \varphi^* f_1(x) + s_I x.
\]
For all $s \in \mcs(\Phi) \cup \{ 0\} $, we have 
\[
\pi_{s} (w) =  \varphi^* \pi_{s} ( f_1(x)) + s_I  \pi_{s} (x) 
=  \varphi^*  \left(  -\frac{s}{\varphi^*} \   x_{s} \right) + s_I x_{s} = (s_I - s) x_{s}.
\]
Choosing $s = s_I$,  we obtain $\pi_{s_I}(w) = 0$.  
Thus, every element $w \in \Phi^{(1)}(A)$ has at least one zero coordinate.

For every $s_I \in \mcs(\Phi) \cup \{ 0\} $, the number of elements  
$x \in M_1$ with $\pi_{s_I}(x)=0$ is 
\[
\prod_{s \in (\mcs(\Phi) \cup \{ 0\}  ) \setminus \{s_I\} } |M'_{s}| = \frac{|M_1|}{|M'_{s_I}|}  
\]
and so the number of elements $x \in M_1$ with $\pi_{s_I}(x)=0$ 
for some $s_I \in \mcs(\Phi) \cup \{ 0\} $ 
is at most $\sum_{s_I \in  \mcs(\Phi) \cup \{ 0\} } |M_1|/|M'_{s_I}|$.
It follows that 
\[
|\Phi^{(1)}(A)| \leq \sum_{s_I \in \mcs(\Phi) \cup \{ 0\} } \frac{|M_1|}{|M'_{s_I}|} 
<  |M_1|\sum_{s_I \in \mcs(\Phi) \cup \{ 0\} } \frac{\varepsilon_1}{ |\mcs(\Phi)| + 1} 
= \varepsilon_1 |M_1|.  
\]
This completes the initial step of the induction.

\section{Proof of Theorem~\ref{HR2:theorem:main}: The inductive step}
  
Assume that $\ell \in \{ 1,2, \ldots, h-1\}$ and that there exist  
a finite $R$-module  $M_{\ell}$ and a function 
$f_{\ell}:M_{\ell}\rightarrow M_{\ell}$ such that the set 
\[
A_{\ell} = A(M_{\ell}, f_{\ell}) 
= \{ f_{\ell}(x_0) : x_0 \in M_{\ell} \} \cup  \{ f_{\ell}(x_0) + x_0 : x_0 \in M_{\ell} \} 
\]
satisfies 
\[
|\Phi^{(\ell)}(A_{\ell} )| < \varepsilon_{\ell}  |M_{\ell}|
\]
where $\Phi^{(\ell)}(A_{\ell})$ is the set 
of all $w_0 \in \Phi(A_{\ell})$ that have a representation of level at most $\ell$.

For every $w_0 \in \Phi(A_{\ell})$, there is an admissible pair 
of functions $(\alpha, \beta)$ of level at most $\ell$ such that 
\[
w_0 = \sum_{x_0 \in M_{\ell}} (\alpha(x_0) f_{\ell}(x_0) + \beta(x_0)(f_{\ell}(x_0)+x_0)).
\]
Because $M_{\ell}$ and $S(\Phi)\cup\{0\}$ are finite sets, there exist only finitely many 
functions from $M_{\ell}$ to $S(\Phi)\cup\{0\}$, 
and only finitely many admissible pairs of functions  from $M_{\ell}$ to $S(\Phi)\cup\{0\}$.  
Let $n$ be the number of   admissible pairs of
functions on $M_{\ell}$ of level exactly $\ell + 1$ with respect to $\Phi$.
We denote these pairs by $(\alpha_i,\beta_i)$ for $i=1, \ldots, n$.

Let $M'_0 = M_{\ell}$ and, for $i=1,\ldots, n$, 
let $M'_i$ be a finite $R$-module with 
\[
|M'_i| > \frac{ n }{ \varepsilon_{\ell+1} - \varepsilon_{\ell}}.  
\]
Let
\[
M_{\ell+1} = M_{\ell} \oplus \bigoplus_{i=1}^{n } M'_i
=  \bigoplus_{i=0}^{n} M'_i.
\]
If the module $M_{\ell}$ depends only on the linear form $\Phi$ and not on $\Upsilon$, then  
the module $M_{\ell+1}$ also depends only on $\Phi$ and not $\Upsilon$.

For
\[
x = (x_0, x_1,\ldots, x_n)  \in M_{\ell+1}
\]
and $i \in \{ 0, 1,\ldots, n\}$, we define the projection $\pi_i:M_{\ell+1} \rightarrow M'_i$ by 
\[
\pi_i(x) = x_i.  
\] 
Define the function
\[
g_0:M_{\ell + 1} \rightarrow M'_0= M_{\ell}
\]
by
\[
g_0(x)  = f_{\ell}(\pi_0(x) ) = f_{\ell}(x_0).
\]
For $i \in \{ 1,\ldots, n\}$, we construct the function  
\[
g_i:M_{\ell + 1} \rightarrow M'_i 
\]
as follows.
If 
\[
\alpha_i(\pi_0(x)) + \beta_i(\pi_0(x)) = 0
\]
then 
\[
g_i(x) = 0.
\]
Recall that $S(\Phi) \subseteq R^{\times}$.  If 
\[
\alpha_i(\pi_0(x)) + \beta_i(\pi_0(x)) \neq 0
\]
then there is a nonempty subset $I_j$ of $\{1,\ldots, h\}$ such that  
\[
\alpha_i(\pi_0(x)) + \beta_i(\pi_0(x)) = s_{I_j} \in R^{\times}.
\]
In this case, let  
\[
g_i(x) = 
- \left( \frac{\beta_i(\pi_0(x))}{ \alpha_i(\pi_0(x)) + \beta_i(\pi_0(x))} \right) x_i.
\]
Define the function $f_{\ell +1}: M_{\ell +1} \rightarrow M_{\ell +1} $ by 
\[
f_{\ell +1}(x) = (g_0(x), g_1(x), \ldots, g_n(x)) = (f_{\ell}(x_0), g_1(x), \ldots, g_n(x))
\]
where $x_0 = \pi_0(x)$.  The diagram
\begin{center}
\begin{tikzcd}
M_{\ell +1} \arrow{r}{f_{\ell +1}}   \arrow{d}[swap]{\pi_0}  & M_{\ell +1}   \arrow{d}{\pi_0} \\
M_{\ell} \arrow{r}{f_{\ell }}  & M_{\ell} 
\end{tikzcd}
\end{center}
commutes.  
As usual, we consider the set  
\begin{align*}
A_{\ell + 1} & = A(M_{\ell +1}, f_{\ell + 1})  \\
& = \{ f_{\ell + 1}(x) : x \in M_{\ell + 1} \} \cup  \{ f_{\ell + 1}(x) + x : x \in M_{\ell + 1} \}.  
\end{align*}

Let 
\[
w \in \Phi^{(\ell+1)}( A_{\ell + 1} ).
\]
There is an admissible pair of functions 
$(\alpha, \beta)$ of level at most $\ell + 1$ such that 
\[
w = \sum_{x \in M_{\ell + 1}} ( \alpha(x) f_{\ell + 1}(x) 
+  \beta(x) (f_{\ell + 1}(x) + x )).
\]
Applying Lemma~\ref{HR2:lemma:commute} with 
$M = M_{\ell+1}$, $M'_0 = M_{\ell}$, $f = f_{\ell + 1}$, and $f_0 = f_{\ell}$, 
we obtain an admissible pair of functions $(\alpha^*,\beta^*)$ 
on $M_{\ell}$ of level $\ell^* \in \{1,2,\ldots, \ell+1\}$ that represents the element 
$\pi_0(w) \in \Phi(A_{\ell})$.
If $\ell^* \leq \ell$, then $\pi_0(w) \in \Phi^{(\ell)}(A_{\ell})$.  
Because the number of elements in $\Phi^{(\ell)}(A_{\ell})$ is less than 
$\varepsilon_{\ell} |M_{\ell}|$, it follows that the number of elements 
$w  \in \Phi^{(\ell +1)}(A_{\ell+1})$ such that $\pi_0(w) \in \Phi^{\ell)}(A_{\ell})$ 
is less than
\[
\varepsilon_{\ell} |M_{\ell}| \prod_{i=1}^n |M'_i| = \varepsilon_{\ell} |M_{\ell+1}|. 
\]

If  the admissible pair $(\alpha^*,\beta^*)$ has level $\ell^* = \ell + 1$, 
then  $(\alpha^*,\beta^*) = (\alpha_i,\beta_i)$ for some $i \in \{1,\ldots, n\}$.  
Moreover, the admissible pair $(\alpha,\beta)$ must also have level $\ell+1$.  
If $\{y_1,\ldots, y_{\ell+1} \}$ is the support of $(\alpha,\beta)$ in $M_{\ell +1}$ 
and if $z_j = \pi_0(y_j)$ for $j \in \{ 1,\ldots, \ell+1\}$,
then $\{z_1,\ldots, z_{\ell+1} \} =  \{\pi_0(y_1),\ldots, \pi_0(y_{\ell+1}) \}  $ 
is the support of $\pi_0(w)$ in $M_{\ell}$.  
For each $z_j \in \{z_1,\ldots, z_{\ell+1} \}$,  there is a unique $n$-tuple 
$(x_{1,j},\ldots, x_{n,j}) \in \bigoplus_{i=1}^n M'_i$ such that 
\[
y_j = (z_j, x_{1,j},\ldots, x_{n,j}).  
\]
Therefore,   
\[
\alpha_i(z_j) = \alpha(y_j)  \qqand \beta_i(z_j) = \beta(y_j)  
\]
and
\[
\alpha_i(z_j) + \beta_i(z_j)  = \alpha(y_j)  + \beta(y_j)  \neq 0
\]
for all $i \in \{1, \ldots, n\}$ and $j \in \{1,\ldots, \ell + 1\}$.  
We have
\[
w = \sum_{j=1}^{\ell + 1} (\alpha(y_j)f_{\ell +1}(y_j)+\beta(y_j)(f_{\ell +1}(y_j)+y_j)
\]
and 
\begin{align*}
\pi_i(w) 
& = \sum_{j=1}^{\ell + 1} ( \alpha(y_j) \pi_i f_{\ell +1}(y_j)  
+  \beta(y_j) \pi_i(f_{\ell +1}(y_j)+y_j))  \\
& = \sum_{j=1}^{\ell + 1} ( \alpha_i(z_j) g_i(y_j) +  \beta_i(z_j) ( g_i(y_j)+ x_{i.j})) \\  
& = \sum_{j=1}^{\ell + 1} (   (\alpha_i(z_j) + \beta_i(z_j)) g_i(y_j)+   \beta_i(z_j) x_{i.j})) \\
& =  \sum_{j=1}^{\ell + 1} (   (\alpha_i(z_j) + \beta_i(z_j)) 
\left( -\frac{\beta_i(z_j)  }{ \alpha_i(z_j) + \beta_i(z_j) }  \right) x_{i,j} +   \beta_i(z_j) x_{i.j})) \\
& = 0.
\end{align*}
The number of elements $w \in \Phi^{\ell+1)}(A_{\ell+1})$ such that $\pi_i(w) = 0$ is
\[
\prod_{\substack{i'=0 \\i' \neq i}}^n |M'_{i'}| = \frac{|M_{\ell+1}|}{|M'_i|} 
<  \frac{  \varepsilon_{\ell+1} - \varepsilon_{\ell}}{ n } |M_{\ell +1}|.
\]
If  $w \in \Phi^{ (\ell+1)}(A_{\ell+1})$, then either 
$\pi_0(w) \in \Phi^{(\ell)}(A_{\ell})$ or $\pi_i(w) = 0$ for some $i \in \{1,2,\ldots, n\}$, 
and so 
\begin{align*}
|\Phi^{(\ell+1)}(M_{\ell +1} | 
& <   \varepsilon_{\ell} |M_{\ell+1}|  
+ \sum_{i=1}^n   \frac{  \varepsilon_{\ell+1} - \varepsilon_{\ell}}{ n } |M_{\ell +1}|  \\
& = \varepsilon_{\ell+1} |M_{\ell +1}|. 
\end{align*}
This completes the induction and the proof of Theorem~\ref{HR2:theorem:main}.

\section{Proof of Theorem~\ref{HR2:theorem:mainZ} }
As usual, if $c$ is a positive integer, then in the ring $R$ we denote $c\cdot 1_R$ by $c$ 
and $-(c \cdot 1_R)$ by $-c$.  Thus, $c$ is a unit in $R$ if $c\cdot 1_R \in R^{\times}$.
For example, 2 is a unit in $R = \Z/3\Z$.

In the statement of Theorem~\ref{HR2:theorem:mainZ}, the sequences of 
nonzero integral coefficients of the linear forms 
$\Upsilon$ and $\Phi$  satisfy $0 \in \mcs(\Upsilon)$ and $0 \notin \mcs(\Phi)$.
There is an infinite set $\mcm$ of positive integers $m$ such that 
$\gcd(s,m) = 1$ for all $s \in \mcs(\Upsilon) \setminus \{ 0\}$ 
and for all $s \in \mcs(\Phi)$.   
It follows that if $m \in \mcm$ and $R = \Z/m\Z$,  
then $\mcs(\Phi) \subseteq R^{\times}$. 
Moreover, every coefficient of $\Upsilon$ is a unit in $R$, and so 
$\{0,u\} \subseteq \mcs(\Upsilon)$ 
for some $u \in R^{\times}$.   
 
The proof of Theorem~\ref{HR2:theorem:mainZ} is essentially the same as the proof 
of Theorem~\ref{HR2:theorem:main}.   
In the initial step of the inductive proof of Theorem~\ref{HR2:theorem:main}, 
which was the case $\ell = 1$, we constructed  
a module $M_1 = \bigoplus_{s \in \mcs(\Phi) \cup \{ 0\} } M'_s$, where the cardinality 
condition~\eqref{HR2:cardinalityM1}  was 
the only constraint on the choice of the finite $R$-modules $M'_s$.   
In the  proof of Theorem~\ref{HR2:theorem:mainZ}, we choose a subset 
\[
\{m'_s: s\in \mcs(\Phi) \cup \{ 0\} \} \subseteq \mcm
\]
whose elements are pairwise relatively prime and satisfy 
\[
m'_s >  \max\left(  \frac{|\mcs(\Phi)| + 1}{\varepsilon_1} ,c \right).
\]
Let 
\[
M_1 =  \bigoplus_{s \in \mcs(\Phi) \cup \{ 0\} } \Z/m'_s\Z \cong \Z/m_1\Z
\]
 where 
 \[
 m_1 = \prod_{s\in \mcs(\Phi) \cup \{ 0\} } m'_s.
 \]
We again use formulae~\eqref{HR2:gs} and~\eqref{HR2:f1} to construct the functions 
$g_s: M_1 \rightarrow \Z/m'_s\Z$ and $f_1:M_1\rightarrow M_1$.

Similarly, in the inductive step, we start with the module $M_{\ell} = \Z/m_{\ell} \Z$ 
and a function $f_{\ell}: M_{\ell}  \rightarrow M_{\ell}$.
Choosing a set $\{m''_1,\ldots, m''_n\}$ 
of pairwise relatively prime integers in \mcm\  
such that $\gcd(m''_i,m_{\ell}) = 1$ for $i \in \{1,\ldots, n\}$, 
we let 
\[
M_{\ell + 1} =  (\Z/m_{\ell}\Z )\oplus  \bigoplus_{i=1}^n \Z/m''_i\Z \cong \Z/m_{\ell + 1}\Z
\]
where 
\[
m_{\ell + 1} = m_{\ell} \prod_{i=1}^n m''_i.
\] 
We complete the proof by constructing functions $g_i$ and $f_{\ell + 1}$ 
exactly as in the proof of Theorem~\ref{HR2:theorem:main}.

\section{Proof of Theorems~\ref{HR2:theorem:mainMany} 
and~\ref{HR2:theorem:mainManyZ} }

To prove Theorem~\ref{HR2:theorem:mainMany}, 
we consider the linear form 
$\chi = \sum_{k=1}^K \Phi_k$.  
Applying Theorem~\ref{HR2:theorem:main}, we obtain a finite $R$-module $M$ 
and a subset $A$ of $M$ such that 
\[
|\chi(A)| < \varepsilon |M|.  
\]
For every $a^* \in A$ and $k \in \{1,\ldots, K \}$, we have 
\[
\Phi_k(A) +  \sum_{\substack{k'=1 \\k' \neq k}}^K  \Phi_{k'}(a^*, \ldots,a^*) 
\subseteq \chi(A)
\]
and so
\[
|\Phi_k(A)| \leq |\chi(A)| < \varepsilon |M|.  
\]
For all $j \in \{1,\ldots, J\}$ we have $\{0,u_j\} \subseteq \mcs(\Upsilon_j)$ 
for some $u_j \in R^{\times}$, and so $\Upsilon_j(A \cup \{ 0\}) = M$.
This completes the proof of Theorem~\ref{HR2:theorem:mainMany}.

Using the same argument, we deduce Theorem~\ref{HR2:theorem:mainManyZ} from 
Theorem~\ref{HR2:theorem:mainZ}.

\section{Open problems}

\bprob
Theorem~\ref{HR2:theorem:main} gives  a sufficient condition 
on a pair of linear forms $\Upsilon$ and $\Phi$ 
to compel  the existence of a module $M$ and subset $A$ of $M$ such that 
$\Upsilon(A \cup \{0\}) = M$ and $|\Phi(A)| < \varepsilon |M|$.  
Is there a corresponding necessary condition?  Is there a necessary and sufficient condition?

We can ask the analogous questions for each of 
Theorems~\ref{HR2:theorem:mainZ}--\ref{HR2:theorem:mainManyZ}.  
\eprob

\bprob
For every prime number $p$ and positive integer $n$, let
$\mathbf{F}_{p^n}$ be the finite field with $p^n$ elements. 
Let $R = \Z$ or $\mathbf{F}_{p^n}$.     
Let $\Upsilon$ be a $g$-ary linear form with nonzero  coefficients in $R$ 
and with $0 \in \mcs(\Upsilon)$.   
Let $\Phi$ be an $h$-ary linear form with nonzero  coefficients in $R$ 
and with $0 \notin \mcs(\Phi)$.  
Let $\varepsilon > 0$.
We say that a finite $R$-module $M$ has property $(\Upsilon, \Phi,\varepsilon)$ 
if $M$ contains a subset $A$ such that 
\beq                \label{HR2:propertyCELF}
\Upsilon(A) = M \qqand |\Phi(A)| < \varepsilon |M|.
\eeq

Let $R = \Z$.  
What is the smallest integer $m = m(\Upsilon, \Phi,\varepsilon)$ 
such that the finite module $\Z/m\Z$ 
has property $(\Upsilon, \Phi,\varepsilon)$?
If $M = \Z/m\Z$ has this property, describe the set 
\[
\{ |A| : A \subseteq \Z/m\Z \text{ and $A$ satisfies~\eqref{HR2:propertyCELF}} \}.
\]
\eprob

\bprob
Do there exist infinitely many prime numbers $p$ such that the finite field 
$\mathbf{F}_p = \Z/p\Z$ has property $(\Upsilon, \Phi,\varepsilon)$?
Does $\mathbf{F}_p$ have property $(\Upsilon, \Phi,\varepsilon)$ 
for all sufficiently large primes $p$?
\eprob

\bprob
Let $p$ be prime number.  Does there exist a positive integer $n$ such that 
the finite field $\mathbf{F}_{p^n}$ has property $(\Upsilon, \Phi,\varepsilon)$? 
Does the finite field $\mathbf{F}_{p^n}$ have property $(\Upsilon, \Phi,\varepsilon)$ for infinitely many $n$, or for all sufficiently large $n$?
\eprob

\bprob
The linear forms $\Upsilon = t_1-t_2$ and $\Phi = 2t_1-t_2$ 
are an interesting special case.   
We have $0 \in \mcs(\Upsilon) = \{-1,0, 1\}$ 
and $0 \notin \mcs(\Phi) = \{ -1,1,2\}$.  
Let  $\varepsilon = 1/2$.  
Compute the smallest integer $m$ such that a subset $A$ of  $M = \Z/m\Z$ 
satisfies~\eqref{HR2:propertyCELF}.    
If $m$ is a positive integer such that there exists $A \subseteq \Z/m\Z$ 
satisfying~\eqref{HR2:propertyCELF}, compute the size of the largest set $A$ 
that satisfies~\eqref{HR2:propertyCELF}.  
\eprob

\bprob
Haight~\cite{haig73} applied his congruence theorem to construct a set $E$ of 
positive real numbers such that $E - E = \R$ but the sumset $hE$ 
has zero Lebesgue measure for all positive integers $h$.  Let $\Phi(t_1,t_2) = 2t_1 - t_2$.  
Does there exist a set $E$ of positive real numbers such that $E - E = \R$ 
but the set $\Phi(E)$ has zero Lebesgue measure?  
\eprob

\bprob
For every finite subset $A$ of an abelian group, we have the 
Freiman-Pigaev~\cite{frei-piga73} inequality  
\[
|A+A|^{3/4} \leq |A-A| \leq |A+A|^{4/3}
\]
Equivalently, with $\Upsilon(t_1,t_2) = t_ 1 - t_2$ and $\Phi(t_1,t_2) = t_ 1 + t_2$, 
we have 
\[
\frac{3}{4} \leq \frac{\log |\Upsilon(A)|}{\log |\Phi(A)|} \leq \frac{4}{3}.  
\]
Let $\Upsilon$ and $\Phi$ be linear forms with integer coefficients such that 
$0 \in S(\Upsilon)$ and $0 \notin S(\Phi)$.  
By Theorem~\ref{HR2:theorem:mainZ}, for every $\varepsilon $ such that 
$0 < \varepsilon < 1$, 
there exist infinitely many positive integers $m$ and subsets $A$ of $\Z/m\Z$ 
such that 
\[
 \frac{\log |\Upsilon(A)|}{\log |\Phi(A)|} > \frac{\log |M|}{\log|M| - \log (1/\epsilon) }.
\]

Can  the Haight-Ruzsa method improve estimates for 
$\log |\Upsilon(A)|/\log |\Phi(A)| $?
\eprob

The next problem is suggested by the following result.  

\bl
Let $M_1, \ldots, M_n$ be finite $R$-modules.  If $M_j$ 
has property $(\Upsilon, \Phi, \varepsilon_j)$ for all $j=1,\ldots, n$, 
then the $R$-module $M_1 \oplus \cdots \oplus M_n$  
has property $(\Upsilon, \Phi, \varepsilon_1\cdots \varepsilon_n)$.
\el

\begin{proof}
For $j =1,\ldots, n$, let $A_j$  be a finite subset of $M_j$ such that 
$\Upsilon(A_j) = M_j$ and  $|\Phi(A_j)| < \varepsilon |M_j|$.  
Let $x = (x_1,\ldots, x_n) \in M_1 \oplus \cdots \oplus M_n$, 
and let   $ (a_{1,j}, \ldots, a_{g,j} ) \in A_j^g$ satisfy  
\[
x_j = \Upsilon(a_{1,j}, \ldots, a_{g,j}) =   \sum_{i=1}^g \upsilon_i a_{i,j} .
\]
Let $A =  A_1 \oplus \cdots \oplus A_n$.
We have 
\begin{align*}
(x_1,\ldots, x_n) 
& = \left(\Upsilon(a_{1,1}, \ldots, a_{g,1}), \ldots, \Upsilon(a_{1,n}, \ldots, a_{g,n})   \right) \\
& = \left(  \sum_{i=1}^g \upsilon_i a_{i,1}, \ldots,   \sum_{i=1}^g \upsilon_i a_{i,n} \right)  \\
& =   \sum_{i=1}^g \upsilon_i  \left( a_{i,1}, \ldots, a_{i,n} \right)  \\
& \in \Upsilon(A) 
\end{align*}
and so $\Upsilon(A) = M_1 \oplus \cdots \oplus M_n$.  

Similarly,  $x = (x_1,\ldots, x_n) \in \Phi(A)$ if and only if  $x_j \in \Phi(A_j)$ for all $j=1,\ldots, n$.
It follows that 
\[
\Phi(A) = \Phi(A_1) \times \cdots \times \Phi(A_n)
\]
and so
\[
|\Phi(A)| = | \Phi(A_1)| \cdots | \Phi(A_n) | 
< \varepsilon_1|A_1| \cdots \varepsilon_n |A_n| 
=  \varepsilon_1 \cdots \varepsilon_n |A|.
\]
This completes the proof.  
\end{proof}

\bprob
Suppose that the finite $R$-modules $M_1$ and $M_2$ have
property $(\Upsilon, \Phi, \varepsilon)$.  
Does there exist $\varepsilon' > 0$ (with $\varepsilon'$ depending only on $\varepsilon$) 
such that the tensor product $M_1 \otimes M_2$ 
has property $(\Upsilon, \Phi, \varepsilon')$?
\eprob

\def\cprime{$'$} \def\cprime{$'$} \def\cprime{$'$} \def\cprime{$'$}
\providecommand{\bysame}{\leavevmode\hbox to3em{\hrulefill}\thinspace}
\providecommand{\MR}{\relax\ifhmode\unskip\space\fi MR }
\providecommand{\MRhref}[2]{%
  \href{http://www.ams.org/mathscinet-getitem?mr=#1}{#2}
}
\providecommand{\href}[2]{#2}

\end{document}